\documentclass[]{amsart}
\usepackage{amscd}
\usepackage{hyperref}
\usepackage[T1]{fontenc}
\hypersetup{
    pdftitle={Surface Homeomorphisms with zero dimensional singular set},
    pdfauthor={Christian Bonatti and Boris Kolev},
    pdfsubject={Mathematics},
    pdfkeywords={Ker\'{e}kj\'{a}rt\'{o} theory; Regular homeomorphisms; Limit set; Riemann sphere}
    }
\newtheorem{thm}{Theorem}[section]
\newtheorem{cor}[thm]{Corollary}
\newtheorem{lem}[thm]{Lemma}

\theoremstyle{definition}

\theoremstyle{remark}
\newtheorem{rem}[thm]{Remark}
\numberwithin{equation}{section}
\newcommand{\norm}[1]{\left\Vert#1\right\Vert}

\newcommand{\set}[1]{\left\{#1\right\}}

\DeclareMathOperator{\cl}{cl\,}
\begin{document}

\title[Surface Homeomorphisms]{Surface Homeomorphisms with zero dimensional singular set}

\author[C. Bonatti]{Christian Bonatti}
\address{Universit\'{e} de Bourgogne, D\'{e}pt. de Math\'{e}matiques, B.P. 138, 21004 Dijon Cedex, FRANCE.}
\email{bonatti@satie.u-bourgogne.fr}%

\author[B. Kolev]{Boris Kolev}%
\address{CMI, 39, rue F. Joliot-Curie, 13453 Marseille cedex 13, France}%
\email{kolev@cmi.univ-mrs.fr}%

\thanks{The authors wish to express their gratitude to John Guaschi, Toby
Hall and Marie-Christine P\'{e}rou\`{e}me for the discussions that helped to
improve this paper, to Lucien Guillou for having introduced them to
the work of Ker\'{e}kj\'{a}rt\'{o} and to the referee of the first version of
this paper for all the precious references he gave to us.}


\date{May 26, 1997}

\begin{abstract}
We prove that if $f$ is an orientation-preserving
homeomorphism of a closed orientable surface $M^{2}$ whose
singular set $\Sigma(f)$ is totally disconnected, then $f$ is
topologically conjugate to a conformal transformation.
\end{abstract}

\maketitle

\section{Introduction}\label{sec:Intro}

In a series of papers \cite{Ker3, Ker5, Ker4, Ker6, Ker8, Ker7},
Ker\'{e}kj\'{a}rt\'{o} gave necessary and sufficient conditions for an
orientation-preserving surface homeomorphism to be conjugate to a
conformal isomorphism, answering a question of Brouwer \cite{Bro1}.
For that purpose, he introduced the notion of \textit{regularity}.

Let $f:X \to X$ be a homeomorphism of a compact metric space. A
point $x\in X$ is \textit{regular} if the family $\set{f^{n} }$ of
all (positive and negative) iterates of $f$ is equicontinuous at
$x$, that is for all $\varepsilon >0$ there exists $\delta >0$ so
that if $d(x,y)<\delta$ then $d(f^{n}(x), f^{n}(y)) < \varepsilon$
for all $n$. This property is obviously independent of the choice
of the metric and invariant under topological conjugacy.

The purpose of this paper is to give a complete exposition of
Ker\'{e}kj\'{a}rt\'{o}'s theory in a more general setting, with shorter and
simpler proofs. Some of Ker\'{e}kj\'{a}rt\'{o}'s arguments have already been
clarified \cite{Bre, Fol, Hem, HL, Lam, Rit}. However, up to the
authors' knowledge, there is nowhere in the literature a complete,
modern and elementary paper on Ker\'{e}kj\'{a}rt\'{o}'s theory. This is what we
intend to do in the following.

Let $f$ be a homeomorphism of a closed surface $M^{2}$. We define
the \textit{singular} set of $f$ , noted $\Sigma(f)$ to be the
closure of the set of non regular points. One may ask what the
connection is between $\Sigma(f)$ and the Julia set $J(f)$ of a
holomorphic map of the Riemann sphere. For us, the only difference
is that a point is regular if the family of all iterates of $f$,
positive and negative, is equicontinuous whereas in the
holomorphic case only positive iterates are required to form an
equicontinuous family. The main result of this paper is the
following:
\begin{thm}\label{thm:mainthm}
Let $f$ be an orientation-preserving homeomorphism of a closed
orientable surface $M^{2}$ whose singular set, $\Sigma (f)$ , is
totally disconnected. Then:
\begin{enumerate}
    \item If $M^{2}=S^{2}$, $f$ is conjugate to a linear fractional transformation.
    \item If $M^{2}= T^{2}$, $f$ is periodic or $f$ is conjugate to the map $(s,t) \mapsto (s + \alpha, t + \beta)$.
    \item If $\chi(M^{2})< 0$, $f$ is periodic.
\end{enumerate}
\end{thm}
Using the fact that any periodic, orientation-preserving
transformation of a closed surface is conjugate to a conformal
transformation \cite{Bro4,Eps}, Theorem \ref{thm:mainthm} can be
restated by saying that an orientation-preserving homeomorphism of a
closed orientable surface is conjugate to a conformal transformation
if and only if $\Sigma (f)$ is totally disconnected.

In order to be complete, we have also treated the case of
orientation-reversing homeomorphisms, of non-orientable closed
surfaces and of surfaces with boundary, orientable or not. We will
show in particular that a homeomorphism of a surface with negative
Euler characteristic is always periodic if its singular set is
totally disconnected. We will give then the complete
\ref{thm:mainthm} classification up to topological conjugacy of non
periodic homeomorphisms with totally disconnected singular set on
surfaces with non negative Euler characteristic: the sphere $S^{2}$,
the projective plane $\mathbb{P}^{2}(\mathbb{R})$, the disc
$\mathbb{D}^{2}$, the torus $\mathbb{T}^{2}$, the Klein bottle
$\mathbf{K}$, the annulus and the M\"{o}bius band. The essential
argument, once we know the result for closed orientable surface, is
to reduce to that case by considering the two-fold orientation
covering or by taking the double of the surface to remove boundary.

In general $\Sigma(f)$ is very large (and so the hypothesis that
$\Sigma(f)$ is totally disconnected is very strong). There are
very simple homeomorphisms of the sphere $S^{2}$ with
$\Sigma(f)=S^{2}$. Let us present such an example.

Let $h$ be the homeomorphism of the plane $\mathbb{R}^{2}$ which
leaves each circle $c_{r}, r \in [0, +\infty[$ centered at $(0, 0)$
globally invariant and whose restriction to $c_{r}$ is the
rotation by angle $\phi(r)$ where $\phi: [0, +\infty[ \rightarrow
\mathbb{R}$ is any continuous function which is constant on no
interval. Let $f$ be the homeomorphism of $S^{2}$ obtained by
extending $f$ to infinity by a fixed point. Then, $f$ is a
homeomorphism of $S^{2}$ such that $\Sigma(f) = S^{2}$.

Another interesting example is the case of a
$C^{1}$-structurally-stable diffeomorphism $f$ of a surface
$M^{2}$ (that is, $f$ satisfies Axiom A and the strong
transversality condition). In that case, a point is regular if and
only if it belongs to the intersection of the stable manifold of a
periodic attractor with the unstable manifold of a periodic
repeller. Indeed, a point which does not lie in such an
intersection must be a limit point of stable or unstable manifolds
of saddle points (according to the Shadowing Lemma, \cite{Newh}) but
these manifolds belong to the singular set (according to the
$\lambda$-Lemma, \cite{Shu}). There are examples of
$C^{1}$-structurally-stable diffeomorphisms $f$ on any surface
$M^{2}$, which do not have any periodic attractor; each attractor
of $f$ is a non-trivial hyperbolic attractor. In that case $\Sigma
(f) = M^{2}$.

However, it is remarkable that such a simple condition ($f$ is a
homeomorphism with $\Sigma(f)$ totally disconnected) implies such
a strong rigidity (f is conjugate to a conformal transformation).
In the same spirit (but in the opposite direction), let us recall
a beautiful result obtained independently by Hiraide \cite{Hir} and
Jorge Lewowicz \cite{Lew}. A homeomorphism is expansive if there is a
constant $\alpha > 0$ such that, for any pair of distinct points
$(x, y)$ there exists $n \in \mathbb{Z}$ such that the distance
between $f^{n}(x)$ and $f^{n}(y)$ is greater than $\alpha$.
Hiraide and Lewowicz showed that any expansive homeomorphism of a
compact surface is conjugate to a pseudo-Anosov homeomorphism.

In the next section, we briefly review some well-known topological
facts. Section~\ref{sec:CompactMetricSpace} contains general results
on the dynamics of a homeomorphism of a compact metric space. In
section~\ref{sec:PositiveEuler}, we study the case of surfaces with
positive Euler characteristic. The main part of this section is
devoted to the sphere. Section~\ref{sec:NegativeEuler} contains the
complete classification of regular homeomorphisms of surfaces with
non positive Euler characteristic and in particular of the torus.


\section{Preliminaries}\label{sec:preliminaries}

Let $( X,d)$ be a compact metric space, and let $2^{X}$ be the set
of all nonempty closed subsets of $X$ . The Hausdorff distance on
$2^{X}$ is defined by
\begin{equation}
    d_{H}(A, B) = \max \set{ \sup_{a\in A} d (a, B), \,  \sup_{b\in B} d (b, A) }.
\end{equation}
With this distance, $2^{X}$ is a compact metric space \cite{Nad}.
Let $(A_{n})$, be a sequence of nonempty closed subsets of $X$. We
define $\liminf (A_{n})$ to be the set of points $x \in X$ such that
every neighborhood of $x$ meets $A_{n}$, for all but a finite number
of values of $n$. In other words, a point $x$ belongs to $\liminf
(A_{n})$ if and only if there exists a sequence $x_{n} \in A_{n}$,
converging to $x$. Similarly, we define $\limsup (A_{n})$ to be the
set of points $x \in X$ such that every neighborhood of $x$ meets
$A_{n}$, for an infinite number of values of $n$. In other words,
$x$ belongs to $\limsup (A_{n})$ if and only if there exists a
sequence $x_{n} \in A_{n}$, admitting a subsequence $x_{n_{i}}$
converging to $x$. These two sets are closed and $\liminf (A_{n})
\subset \limsup (A_{n})$. The sequence $(A_{n})$ is convergent for
the Hausdorff metric if and only if $\liminf (A_{n}) = \limsup
(A_{n})$ \cite{Nad}.

A continuum is a nonempty, compact, connected metric space. It is
non-degenerate provided that it contains more than one point. The
proof of the following Lemma may be found in \cite{Nad}.
\begin{lem}\label{lem:SpaceContinuaClose}
Let $(A_{n})$ be a sequence of continua of X such that $\liminf
(A_{n})\neq \emptyset$, then $\limsup (A_{n})$ is a continuum. In
particular the space of continua of $X$, denoted $C(X)$, is closed
in $2^{X}$.
\end{lem}
Let $f$ be a homeomorphism of a compact metric space ($X,d)$.
Recall that a point $x \in X$ is regular if the family ${f^{n}}$
is equicontinuous at $x$. For every $\varepsilon > 0$, we let
$\varphi(x, \varepsilon)$ be the least upper bound of positive
numbers $\delta$ such that
\begin{equation}
d(x,y) < \delta \Rightarrow d(f^{n}(x),f^{n}(y)) < \varepsilon,
\quad \forall n.
\end{equation}
By its very definition $\varphi(x, \varepsilon)\leq \varepsilon$.
In the introduction, we have defined the singular set of $f$ ,
written $\Sigma(f)$ to be the closure of the set of points which
are not regular. This set is clearly invariant under $f$ and
$\Sigma(f^{n}) = \Sigma (f)$ for all $n \neq 0$.

If $\Sigma(f)$ is empty, the family $\set{ f^{n}, n \in \mathbb{Z}
}$ is uniformly equicontinuous on $X$. Such homeomorphisms are
called \textit{regular} by Ker\'{e}kj\'{a}rt\'{o} and they are described by the
following well-known result \cite{GH,WD}.

\begin{thm}\label{thm:Equicontinuity}
Let $f$ be a homeomorphism of a compact metric space $(X ,d)$. The
following statements are equivalent:
\begin{enumerate}
    \item The family $\set{ f^{n}, n \in \mathbb{Z} }$ is equicontinuous on $X$.
    \item The closure $G$ of the family $\set{f^{n}, n \in \mathbb{Z} }$ is a compact subgroup of $Homeo(X)$
    with the $C^{0}$-uniform topology.
    \item There exists a metric on $X$ compatible with its topology in which $f$ is an isometry.
\end{enumerate}
\end{thm}

\begin{rem}
$1) \Rightarrow 2)$ is essentially Ascoli's
theorem, whereas in $2) \Rightarrow 3)$ an invariant metric can be
constructed using Haar's measure on $G$.
\end{rem}

A homeomorphism of a compact metric space is \textit{recurrent}
provided that for all $\varepsilon > 0$, there exists $n \neq 0$
such that $d(f^{n}, Id) < \varepsilon$. Clearly, in view of
Theorem~\ref{thm:Equicontinuity}, a regular homeomorphism of a
compact metric space is recurrent.


\section{Dynamics on a compact metric
space}\label{sec:CompactMetricSpace}

In this section, we investigate the connection between regularity
and the limit sets of a homeomorphism $f$ of a compact metric space
$X$ . Our primary interest is in such homeomorphisms which have a
totally disconnected singular set. Recall that a space $X$ is
totally disconnected, provided any connected subset of it is a point
or empty. The subject of this section may be compared to the one of
\cite{HL, Lam}.

Given $x \in X$, we set $\mathcal{O}(x, f) = \set{f^{n}(x) ;\, n \in
\mathbb{Z} }$. We define similarly $\mathcal{O}^{+} (x, f)$ as the
set of positive iterates of $x$ and $\mathcal{O}^{-}(x,f)$ as the
set of negative iterates of $x$. The $\omega$-limit set of a point
$x \in X$ is defined by $\omega (x, f) = \set{\lim f^{n_{k}}(x),
n_{k} \rightarrow +\infty }$ and its $\alpha$-limit set by a
$\alpha(x, f ) = \omega (x, f^{-1})$. The limit set of $x$ is
$\lambda (x, f) = \alpha (x, f) \cup \,\omega (x, f)$. These sets
are invariant under $f$. A point $x$ is recurrent if $x \in \lambda
(x, f)$.

\begin{lem}\label{lem:FourProperties}
We have the four following properties:
\begin{enumerate}
    \item Let $x \in X$ be a regular point and suppose that there are sequences
    $(x_{i})$ and $(n_{i})$ such that $x_{i}\rightarrow x$
    and $f^{n_{i}}(xi) \rightarrow z$, then $f^{n_{i}}(x) \rightarrow z$.
    \item Let $x$ and $y$ be two regular points. Then $f^{n_{i}}(x)\rightarrow y$ iff $f^{-n_{i}} (y) \rightarrow x$.
    \item A regular point $x$ is recurrent iff $\lambda(x,f)$ contains a regular point.
    \item If a regular point is recurrent, then $\omega(x, f)= \alpha(x, f)= \cl (\mathcal{O}(x, f))$.
\end{enumerate}
\end{lem}

\begin{proof}
1) Let $\varepsilon> 0$ and $\delta = \varphi(x,
\varepsilon)$. For $i$ large enough, we have $d(x_{i}, x) <
\delta$ and $d(f^{n_{i}} (x_{i}), z) < \delta$ and hence
$d(f^{n_{i}}(x), z ) < 2\epsilon$.

2) Let $\varepsilon> 0$ and $\delta = \varphi (x, \varepsilon)$.
For $i$ large enough, we have $d(f^{n_{i}} (x), y) < \delta$ and
hence $d(x, f^{-n_{i}}(y)) < \varepsilon)$.

3) Let $y \in \lambda (x,f)$ be a regular point. Let $\varepsilon>
0$ and $\delta = \varphi(y, \varepsilon )$ . We can find $n> m >
0$ such that $d(y, f^{n}(x)) < \delta$ and $d(y, f^{m}(x)) <
\delta$. Hence, $d(f^{-m}(y), f^{n-m}(x)) < \varepsilon$ and
$d(f^{-m}(y), x)<\varepsilon$ and so $d(x, f^{n-m}(x))
<2\varepsilon$.

4) If $x$ is recurrent then $x \in \omega(x, f)$ and $x \in \alpha
(x, f )$ according to 2). Hence $\cl O(x, f) =\omega(x,f)
=\alpha(x,f)$.
\end{proof}

\begin{cor}\label{cor:RecurrentOpenClose}
The set of points of $X \setminus \Sigma(f)$ which are recurrent is
open and closed in $X \setminus \Sigma(f)$.
\end{cor}

\begin{proof}
According to Lemma~\ref{lem:FourProperties}(1), the set of recurrent
points of $X \ \Sigma(f)$ is closed in $X \ \Sigma(f)$. Let $x \in X
\ \Sigma(f)$ be a recurrent point. We can find a neighborhood $U
\subset X \setminus \Sigma(f)$ of $x$ such that $\lambda(y, f) \cap
X \setminus \Sigma(f) \neq\emptyset$ for each $y \in U$. According
to Lemma~\ref{lem:FourProperties}(3), each point of $U$ is recurrent
so the set of recurrent points of $X \setminus \Sigma(f )$ is also
open.
\end{proof}

\begin{lem}\label{lem:RecurrentOrbit}
Let $x \in X \ \Sigma(f)$. If $x$ is recurrent, then for every
neighborhood $U$ of $x$ there exists an integer $N \geq 0$ such that
\begin{equation}
  \mathcal{O}(x, f) \subset \bigcup_{i=0}^{N} f^{i}(U).
\end{equation}
\end{lem}

\begin{proof}
Let $x \in X \setminus \Sigma(f)$ be a recurrent point. According to
Lemma~\ref{lem:FourProperties}(4), we have
\begin{equation}
\cl(\mathcal{O}(x,f)) = \omega(x,f) =\alpha(x,f).
\end{equation}
Let $U$ be a compact neighborhood of $x$ such that $U \subset X
\setminus \Sigma (f)$ and let $V \subset U$ be an open neighborhood
of $x$. For each point $y \in \lambda (x, f )$, we can find by
Lemma~\ref{lem:FourProperties}(2) an integer $n(y)> 0$ such that
$f^{n(y)}(y) \in V$ and hence we can find an open neighborhood
$V_{y}$ of $y$ such that
\begin{equation}
f^{n(y)} (V_{y}) \subset V.
\end{equation}
Let $V_{1}, V_{2}, \dotsc, V_{r}$ be a finite subcover of the
covering $(V_{i})$ of $\lambda(x, f) \cap U$. To each open set
$V_{i}$ there corresponds a positive integer $n_{i}$ such that
\begin{equation}
f^{n_{i}}(V_{i}) \subset V
\end{equation}
We let $N = \max{n_{i}}$. For any $z \in \lambda(x,f ) \cap U$,
there exists $n \in \set{0, 1, \dotsc, N }$ such that
\begin{equation}
f^{n}(z) \in V.
\end{equation}
Therefore, we can construct an increasing sequence
$m_{i}\rightarrow +\infty$ such that:
\begin{enumerate}
    \item $m_{i+l} - m_{i} \leq N$, for all $i$,
    \item $f^{m_{i}}(x) \in V \subset U$, for all $i$
\end{enumerate}
and we have thus:
\begin{equation}
\mathcal{O}(x, f) \subset \omega(x, f) = \cl (\mathcal{O}^{+} (x,
f)) \subset \bigcup_{i=0}^{N} f^{i} (U).
\end{equation}
\end{proof}

\begin{cor}\label{cor:RegularLimitSet}
Let $x \in X \setminus\Sigma(f)$. Then $\lambda(x, f)\cap \Sigma(f)
\neq \emptyset$ iff $\lambda(x, f) \subset \Sigma(f)$.
\end{cor}

\begin{proof}
Let $x \in X \setminus\Sigma(f)$ and suppose $\lambda (x, f)
\not\subset \Sigma(f)$. According to
Lemma~\ref{lem:FourProperties}(3), $x$ is recurrent and from
Lemma~\ref{lem:RecurrentOrbit}, we get that $\lambda (x,f ) = \cl
(\mathcal{O}(x, f)) \subset X \setminus \Sigma (f )$.
\end{proof}

\begin{thm}\label{thm:SingularSet}
Let $f$ be a homeomorphism of a locally connected, compact metric
space $( X ,d )$ such that $\Sigma(f)$ is totally disconnected.
Then:
\begin{enumerate}
    \item For each $s \in \Sigma(f)$ which is not a regular point, there exists $x
    \in X \setminus \Sigma (f)$ such that $s$ belongs to the limit set
    $\lambda (x, f)$ of $x$.
    \item For each $x_{0} \in X \setminus \Sigma (f)$ and any sequence $n_{k} \in \mathbb{Z}$ such that
$\lim f^{n_{k}} (x_{0}) = s \in \Sigma(f)$, one has in fact
$f^{n_{k}}(x) \rightarrow s$ uniformly on a neighborhood $U$ of
$x_{0}$ and the set
\begin{equation}
E_s  = \set{ {x \in X\setminus \Sigma \left( f \right);\quad \lim
f^{n_k } \left( x \right) = s} }
\end{equation}
is open and closed in $X\setminus\Sigma (f)$
\end{enumerate}
\end{thm}

\begin{proof}
1) Suppose that $s \in \Sigma(f)$ is not a regular
point. There exist $\varepsilon > 0$, a sequence $z_{p}
\rightarrow s$ and a sequence $n_{p}$ of integers so that
 \begin{equation}
d\left(f^{n_{p}}(z_{p}), f^{n_{p}}(s)\right) \geq
\varepsilon,\quad \forall p.
\end{equation}
Since $\Sigma(f)$ is a compact, totally disconnected, proper subset
of $X$ (otherwise $X = \Sigma(f)$ is finite and therefore $X =
\Sigma(f) = \emptyset$), we can find a neighborhood $U$ of $\Sigma(f
)$, $U \neq X$, which is the union of non intersecting open sets of
diameter less than $\varepsilon$ (cf 7.10 of \cite{Nad}). Let
\begin{equation}
K = \cl\left(\bigcup_{n\in\mathbb{Z}} f^{n} (X \setminus
U)\right).
\end{equation}

Suppose that $s \notin K$. As $X$ is locally connected, the
components of $X \setminus K$ are open and since $X \setminus K
\subset U$, all of them are of diameter less than $\varepsilon$.
Furthermore, these components are permuted by $f$ and so for $p$
large enough, $f^{n}(z_{p})$ and $f^{n}(s)$ are in the same
component of $X \setminus K$ which implies $d(f^{n} (z_{p}),f^{n}
(s)) < \varepsilon$ for all $n$. This is a contradiction.

Hence $s \in K$ and we can find sequences $x_{n} \in X \setminus U$
and in $i_{n}\in \mathbb{Z}$ such that $f^{i_{n}} (x_{n})$ converges
to $s$. Extracting a subsequence if necessary we can assume that
$x_{n}$ converges to a point $x \in X \setminus U$ as $n \rightarrow
+\infty$. From Lemma~\ref{lem:FourProperties}, one obtains that
$f^{i_{n}}(x) \rightarrow s$ as $n \rightarrow +\infty$, that is $s
\in \lambda(x, f )$.

2) Let $x_{0} \in X \setminus \Sigma (f)$, and suppose that
$f^{n_{k}}(x_{0}) \rightarrow s \in \Sigma(f)$ as $k \rightarrow
+\infty$. According to Corollary~\ref{cor:RegularLimitSet}, $x_{0}
\notin \omega(x_{0}, f)$. Choose a compact and connected
neighborhood $U$ of $x_{0}$ such that $x \notin \omega(x, f)$ for
all $x \in U$ (cf Corollary~\ref{cor:RecurrentOpenClose} and
Lemma~\ref{lem:FourProperties}(4)). According to
Lemma~\ref{lem:SpaceContinuaClose} $\limsup (f^{n_{k}}(U))$ is a
continuum since $s \in \liminf (f^{n_{k}}( U ) )$. If $\limsup
(f^{n_{k}}(U))$ contains some point $y \in X \setminus\Sigma(f )$,
there exists a subsequence $m_{k}$ of $n_{k}$, and a sequence $x_{k}
\in U$ such that $f^{m_{k}} (x_{k})$ converges to $y$. From
compactness of $U$, if necessary extracting a subsequence, we may
assume that $x_{k}$ converges to a point $x \in U$, and from
Lemma~\ref{lem:FourProperties}(1), we have $\lim f^{m_{k}} (x) = y$
which leads to a contradiction since $\omega (x, f) \subset \Sigma
(f)$ for all $x \in U$. Hence
\begin{equation}
    \limsup(f^{n_{k}}(U)) \subset \Sigma(f)
\end{equation}
and is therefore reduced to $\set{s }$ since $\Sigma(f)$ is totally
disconnected. This shows in particular that the set $E_{s}$ is open
in $X \setminus \Sigma(f )$. That it is also closed results from
Lemma~\ref{lem:FourProperties}(1).
\end{proof}

\begin{cor}\label{cor:SingularSetEmpty}
Let $f$ be a homeomorphism of a closed surface $M^{2}$ with a
totally disconnected singular set $\Sigma(f )$. Then $\Sigma(f)$ is
empty unless $M^{2} = S^{2}$ and in that case $\Sigma(f)$ contains
no more than two points.
\end{cor}

\begin{proof}
Suppose that $\Sigma(f )$ is not empty and let $s \in \Sigma(f)$ be
a non regular point. According to Theorem~\ref{thm:SingularSet},
there exists $x_{0} \in M^{2} \setminus \Sigma(f)$ such that $s \in
\lambda(x_{0}, f)$. To fix ideas, let us suppose that $s \in \omega
(x_{0}, f)$ and we choose a sequence $n_{k} \rightarrow +\infty$ so
that $\lim f^{n_{k}} (x_{0}) = s$. According to
Theorem~\ref{thm:SingularSet}, the set
\begin{equation}
E_{s} = \set{x \in M^{2}\setminus \Sigma(f) ; \quad \lim
f^{n_{k}}(z) = s }
\end{equation}
is open and closed in $M^{2} \setminus\Sigma(f)$ and the
convergence is uniform on every compact subset of $E_{s}$. Since
the connectedness of a surface is preserved by the removing of a
totally disconnected subset, $E_{s} = M^{2} \setminus \Sigma(f )$.
Therefore, $\set{ s  } \subset \omega(x, f) \subset \Sigma(f)$ for
all $x \in M^{2}\setminus \Sigma(f )$. Let $\alpha$ be an arc
joining $x$ and $f(x)$ in $M^{2} \setminus \Sigma(f )$, then:
\begin{equation}
    \limsup(f^{n}(\alpha)) = \bigcap_{n=0}^{+\infty}\overline{\bigcup_{k\geq n} f^{k}(\alpha)}
\end{equation}
is a continuum which lies entirely in $\Sigma(f)$. Otherwise, there
would be a point $y \in \alpha$ such that $\omega(y, f)\not\subset
\Sigma(f)$ which is impossible by
Corollary~\ref{cor:RegularLimitSet}. Hence, $\limsup(f^{n}(\alpha))
= \set{ s  }$ and $\omega(x, f) = \set{s }$ for every $x \in M^{2}
\setminus \Sigma(f)$.

If $\Sigma(f)$ is not reduced to $s$, we can find a non regular
point $s'\neq s$ and we prove similarly that $\alpha (x, f) =
\set{ s'  }$ for all $x \in M^{2} \setminus \Sigma (f )$.
Therefore $\Sigma(f)$ cannot contain more than two points.

The argument above shows that $f^{n}(x) \rightarrow s$ as $n
\rightarrow +\infty$ uniformly on all compact subsets of $M^{2}
\setminus \Sigma (f)$. If $M^{2} \neq S^{2}$, there exists a closed
curve $J$ non homotopic to zero and we can choose such a curve in
$M^{2} \setminus \Sigma(f)$. Let $U$ be a simply connected
neighborhood of $s$. For $n$ large enough, $f^{n}(J) \subset U$ and
hence $f^{n}(J)$ is null homotopic which gives a contradiction.
\end{proof}

\begin{rem}
Corollary~\ref{cor:SingularSetEmpty} can be generalized easily to
any compact manifold $M$ of dimension $dim(M) \geq 2$:
\begin{enumerate}
    \item A compact and totally disconnected subset of a closed manifold
$M$ of dimension $dim(M) \geq 2$ does not separate $M$ . Hence, as
in the proof of Corollary~\ref{cor:SingularSetEmpty}, the singular
set, $\Sigma(f)$, of any homeomorphism of $M$ contains less than two
points or is not totally disconnected.
    \item Let $f$ be a homeomorphism of a closed manifold $M$ of
dimension $dim(M) \geq 2$. If there exists $i \in \set{1,2,\dotsc, n
- 1 }$ such that the homotopy group $\Pi_{i}(M)$ is not $\set{ 0 }$,
and if $\Sigma(f)$ is totally disconnected, then $\Sigma(f) =
\emptyset$. The proof is identical to the proof of
Corollary~\ref{cor:SingularSetEmpty}.
\end{enumerate}
\end{rem}


\section{The case of surfaces with $\chi(M^{2})>
0$}\label{sec:PositiveEuler}

Recall that every conformal automorphism of the Riemann Sphere $\widehat{\mathbb{C}}$ can be expressed as a fractional linear transformation
\begin{equation}
f(z) = \frac{(az + b)}{(cz + d)},
\end{equation}
where the coefficients are complex numbers with determinant $ad - bc \neq 0$. Every non identity automorphism of this type has two distinct fixed points or one double fixed point in $\widehat{\mathbb{C}}$. The non identity automorphisms of $\widehat{\mathbb{C}}$ fall into three classes, as follows:
\begin{itemize}
    \item An automorphism $f$ is said to be \emph{elliptic} if it has two distinct fixed points at which the modulus of the derivative is 1.
    \item The automorphism $f$ is said to be \emph{hyperbolic} if it has two distinct fixed points at which the modulus of the derivative is not 1.
    \item The automorphism $f$ is \emph{parabolic} if it has just a double fixed point.
\end{itemize}
It is an exercise to show that up to topological conjugacy there is only one model of parabolic automorphism, namely the translation $T(z) = z + 1$. There is only one model of hyperbolic automorphism, namely $H(z) = 2z$. But there is a one parameter family of elliptic transformations which are not topologically equivalent, namely the family
$R_{\alpha}(z) = e^{i\alpha z}$.
In the the first three paragraphs of this section we study the case of orientation-preserving homeomorphisms of the sphere and we prove the following theorem.
\begin{thm}\label{thm:MainSphere}
An orientation-preserving homeomorphism of the sphere $S^{2}$ is
topologically conjugate to a linear fractional transformation iff
its singular set contains no non-degenerate continuum or
equivalently if its singular set contains at most two points. More
precisely, the transformation is conjugate to an elliptic, parabolic
or hyperbolic transformation according to whether the number of
singular points is zero, 1 or 2.
\end{thm}
According to Corollary~\ref{cor:SingularSetEmpty}, the singular set
$\Sigma(f)$ of a homeomorphism of the sphere which is totally
disconnected contains no more than two points. The proof of
Theorem~\ref{thm:MainSphere} will be divided into three parts
according as the number of singular points is zero, 1 or 2.


\subsection{The parabolic case}\label{subsec:Parabolic}

\begin{thm}
An orientation-preserving homeomorphism of the sphere with exactly one singular point is
conjugate to the map $z \mapsto z + 1$ of the Riemann sphere.
\end{thm}

\begin{proof}
According to Theorem~\ref{thm:SingularSet}, we have:
\begin{equation}
    \alpha (x,f) = \omega (x,f) = \set{ N  },\quad \forall x\in S^{2}
\end{equation}
and $f^{n}(x) \rightarrow N$ as $n \rightarrow \pm \infty$,
uniformly on every compact subset of $S^{2} \setminus \set{ N  }$.
Hence the group $< f >$ generated by $f$ acts freely and properly
on the plane $\Gamma = S^{2} \setminus \set{ N  }$. The quotient
space $\Gamma_{f}$ is therefore a topological surface and the
projection $\pi_{f}: \Gamma \rightarrow \Gamma_{f}$ is a regular
covering. Since $\Gamma$ is simply connected, the fundamental
group of $\Gamma_{f}$ is isomorphic to $\mathbb{Z}$. Hence,
$\Gamma_{f}$ is homeomorphic to the cylinder $\mathbb{R} \times
S^{1}$. By uniqueness of the universal cover up to isomorphism,
there exist a pair of homeomorphisms $(H,h )$ and a commutative
diagram
\begin{equation}
\begin{CD}
\Gamma @>H>> \mathbb{R}^{2} \\
@VVV        @VVV\\
\Gamma_{f} @>h>> \mathbb{R} \times S^{1}
\end{CD}
\end{equation}

In particular, $H$ is a topological conjugacy between $f_{|\Gamma}$ and the translation $(x, y) \mapsto (x + 1,y)$ of the plane $\mathbb{R}^{2}$.
\end{proof}


\subsection{The hyperbolic case}\label{subsec:Hyperbolic}

\begin{thm}
An orientation-preserving homeomorphism of the sphere with exactly
two singular points is conjugate to the map $z \mapsto 2z$ of the
Riemann sphere.
\end{thm}

\begin{proof}
Let $N$ and $S$ be the two singular points. According to the proof
of Corollary 3.6, $f^{n}(x) \rightarrow N$ for $n \rightarrow
+\infty$ (resp. $f^{n}(x) \rightarrow S$ for $n \rightarrow
-\infty$) uniformly on every compact subset of $S^{2} \ \set{ S }$
(resp. $S^{2} \ \set{ N  }$ ) . This shows that $f$ acts freely
and properly on the cylinder $\Gamma = S^{2} \setminus \set{ N , S
}$. Hence the quotient space $\Gamma_{f}  = \Gamma / f$ is a
topological surface and the projection $\pi_{f}$ is a regular
covering. Since $f$ is orientation-preserving, $\Gamma_{f}$ is
orientable. Let $\gamma \subset \Gamma$ be a simple closed curve
separating S and N . Since $f^{n}(\gamma)$ converges uniformly to
$S$ for positive integers, there exists $n\in \mathbb{N}$ such
that $f^{n}(\gamma)\cap \gamma = \emptyset$. Since $\gamma$ and
$f^{n}(\gamma)$ are two non-intersecting essential simple closed
curves, they bound a closed annulus A on the cylinder $\Gamma$.
Hence $\Gamma_{f^{n}}$ is a torus and since $\Gamma_{f^{n}}$ is a
cyclic regular covering of $\Gamma_{f}$, we also obtain that
$\Gamma_{f}$ is a torus $T^{2}$. Since, up to isomorphism, there
is only one regular covering $\pi : S^{1} \times \mathbb{R}
\rightarrow T^{2}$ whose automorphism group is isomorphic to
$\mathbb{Z}$, there exist a pair of homeomorphisms $(H,h)$ and a
commutative diagram
\begin{equation}
\begin{CD}
\Gamma @>H>> S^{1} \times \mathbb{R} \\
@VVV        @VVV\\
\Gamma_{f} @>h>> T^{2}
\end{CD}
\end{equation}
In particular, $H$ is a topological conjugacy between
$f_{|\Gamma}$ and the translation $(x,y)\mapsto(x, y + 1)$ of the
cylinder $S^{1}\times\mathbb{R}$.
\end{proof}


\subsection{The elliptic case}\label{subsec:Elliptic}

\begin{thm}\label{thm:RegularElliptic}
A regular, orientation-preserving homeomorphism of the sphere is
conjugate to a rotation $R_{\alpha}(z) = e^{i\alpha z}$.
\end{thm}
The first step to this end is to show that around a fixed point of a
regular homeomorphism of the sphere there are arbitrary small
invariant simple closed curves. Compare the results of this
paragraph to those of \cite{Bre, Fol, Hem,Ker5,Rit, Sin}. The
following theorem is a very classical result in elementary conformal
representation theory (see \cite[Theorem 2.6]{Pom}), although it can
be proved also by methods of plane topology \cite{Newm,Why}. Recall
that a point $a \in K$ is a cut point of $K$ if $K \setminus \set{a
}$ is not connected.

\begin{thm}\label{thm:Caratheodory}
Let K be a non-degenerate, locally connected continuum of $S^{2}$
with no cut points. Then the boundary of each component of $S^{2}
\setminus K$ is a simple closed curve.
\end{thm}

\begin{lem}\label{lem:LocalConnexity}
Let $f$ be a regular homeomorphism of the sphere and let $D \subset
S^{2}$ be a closed disc. Then the compact set $K = \cl
(\bigcup_{n\in\mathbb{Z}}f^{n}(D))$ is locally connected.
\end{lem}

\begin{proof}
Let us first recall that a compact metric space is locally connected if and only if for every $\varepsilon > 0$ it is the union of a finite number of compact connected sets of diameter less than $\varepsilon$.

Let $\varepsilon > 0$ be given. Choose a triangulation of $D$ into
a finite number of closed 2-cells, $e_{1},e_{2},\dotsc , e_{r}$,
each of which of diameter less than $\varphi (\varepsilon)$, where
$\varphi (\varepsilon)$ was defined in section 2. For each $n \in
\mathbb{Z}$ we have
\begin{equation}
    f^{n}(D)=\bigcup_{i=1}^{r}e_{i}^{n}
\end{equation}
where $diam (e^{n}_{i}) < \varepsilon$ for $i = 1,\dotsc, r$
($e^{n}_{i}i$ stands for $f^{n}(e_{i})$).

Let $\rho > 0$ so that each 2-cell $e_{i}$ contains a disc $B (x_{i}, \rho)$ in its interior. Then
\begin{equation}
B(f^{n}(x_{i}),\varphi (\rho )) \subset f^{n}(e_{i}),\quad \forall i,\quad \forall n.
\end{equation}
Therefore the family $\set{ e_{i}^{n} }$ contains only finitely
many pairwise non intersecting 2-cells. Let
$\set{e_{i_{1}}^{n^{1}}, \dotsc, e_{i_{p}}^{n_{p}} }$ be a finite
collection of pairwise nonintersecting 2-cells of maximal
cardinality. Then for every $n \in \mathbb{Z}$ and every $i \in
\set{1, \dotsc, r  }$ , there exists $j \in \set{1, \dotsc, p  }$
so that $e_{i_{j}}^{n_{j}}\cap e_{i}^{n}\neq\emptyset$. For each
$k \in \set{1, \dotsc, p  }$, let $M_{k}$ be the closure of the
union of all 2-cells $e_{i_{n}}$ which meet $e_{i_{k}}^{n_{k}}$.
Then $M_{k}$ is a compact connected set of diameter less than
$3\varepsilon$ and
\begin{equation}
    K=\bigcup_{k=1}^{p}M_{k}.
\end{equation}
\end{proof}

\begin{cor}\label{cor:FixedPointNeighborhood}
Let $f$ be a regular homeomorphism of the sphere and let $x$ be a
fixed point for $f$. There exist arbitrarily small closed discs,
invariant under $f$, which are neighborhoods of $x$.
\end{cor}

\begin{proof}
Let $\varepsilon > 0$, $\delta = \varphi (\varepsilon)$
and $\eta =
\varphi(\delta)$. Let $D^{0}$ be the open disc of center $x$ and
radius $\eta$ and set
\begin{equation}
U = \bigcup_{n\in\mathbb{Z}} f^{n} (D^{0}).
\end{equation}
Then $U$ is an open connected subset of $B(x, \delta)$ and $f (U) =
U$. By Lemma~\ref{lem:LocalConnexity}, $\overline{U}$ is locally
connected and has no cut point since the closure of any connected
open set of the sphere has no cut point. According to
Theorem~\ref{thm:Caratheodory}, the boundary of each component of
$S^{2}\setminus \overline{U}$ is a simple closed curve. Since
\begin{equation}
f(S^{2} \setminus B (x, \varepsilon)) \subset S^{2} \setminus B (x, \delta),
\end{equation}
the component of $S^{2} \setminus \overline{U}$ which contains
$S^{2} \setminus B(x,\delta)$ is invariant under $f$. Its boundary
$\gamma$ is an invariant simple closed curve and the open disc
$\Delta^{0}$ bounded by $\gamma$ and containing $x$ is an invariant
disc contained in $B(x,\varepsilon)$.
\end{proof}

An orientation-preserving homeomorphism $f$ of $S^{2}$ has at least
one fixed point. If furthermore, $f$ is regular, there exists an
invariant disc around this point. According to Brouwer's fixed point
theorem, there is another fixed point and up to conjugacy, we can
suppose that $f$ fixes the two poles $N$ and $S$ of the sphere
$S^{2}$.

Let $\gamma$ be an invariant simple closed curve under $f$. We will
denote the rotation number of the restriction of $f$ to $\gamma$ by
$\rho(\gamma, f)$. Readers not familiar with rotation numbers may
refer to \cite{Dev} for an excellent exposition of this notion.

\begin{lem}\label{lem:RotNumberIdentity}
Let $f$ be a regular homeomorphism of the sphere and suppose that
there exists an invariant simple closed curve $\gamma$ such that
$\rho(\gamma, f) = 0$. Then $f$ is the identity map of $S^{2}$.
\end{lem}

\begin{proof}
Recall that if $\rho(\gamma, f)=0$, $f$ has a fixed point on
$\gamma$. The argument below in one dimension lower shows that a
regular orientation-preserving homeomorphism of the circle with a
fixed point is the identity map. Hence, $\gamma \subset Fix(f)$. We
will now prove that the connected component $\Gamma$ of $\gamma$ in
$Fix(f)$ is open and closed in $S^{2}$ which will complete the
proof.

Let $x \in \Gamma, \varepsilon > 0 , (2\varepsilon < diam (\gamma)$)
and $\delta = \varphi (\varepsilon )$. For each ball $B (x, r)$ with
$r < \delta$, we construct as in the proof of
Corollary~\ref{cor:FixedPointNeighborhood} an invariant simple
closed curve $\gamma_{r}$ which bounds a disc $\Delta_{r}\subset
B(x, \varepsilon)$ containing $x$. Since $\gamma_{r}$ meets
$\Gamma$, $f_{|\gamma_{r}} = Id$ and hence $\gamma_{r}\subset
\Gamma$. But
\begin{equation}
\gamma_{r} \subset \cl ( \bigcup_{n\in\mathbb{Z}} f^{n}(c_{r}))
\end{equation}
where $c_{r}$ is the boundary of $B (x, r)$. Hence, for every point
$y \in\gamma_{r}$, there exists a point $x \in c_{r}$ and a sequence
$(n_{k})$ so that $y = lim f^{n_{k}}(x)$ and we have $x = lim
f^{-n_{k}} (y) = y$ (cf. Lemma~\ref{lem:FourProperties}). Therefore
$\gamma_{r} =c_{r}$ and $c_{r}\subset\Gamma$. Since $r$ is
arbitrary, this shows that $B (x, \delta) \subset \Gamma$.
\end{proof}

We now fix an invariant simple closed curve $\gamma$ and let $G$ be
the closure of the group generated by $f$ in $Homeo(S^{2})$. As we
have seen in Theorem~\ref{thm:Equicontinuity}, $G$ is a compact
commutative subgroup of regular homeomorphisms. Each element of $G$
leaves the curve $\gamma$ invariant and moreover we have:

\begin{lem}\label{lem:RotNumberMonomorphism}
The map $\rho_{\gamma}: G \rightarrow S^{1}$ given by $g\mapsto
\rho(\gamma, g)$ induces a bi-continuous isomorphism between $G$ and
a compact subgroup of $S^{1} = \mathbb{R}/\mathbb{Z}$.
\end{lem}

\begin{proof}
It is a classical fact that the rotation number of a circle map
depends continuously on the map. Therefore $\rho_{\gamma} : G
\rightarrow S^{1}$ is continuous and from the relation
\begin{equation}
\rho_{\gamma}(f^{n}) =n\rho_{\gamma}(f),\quad \forall n \in\mathbb{Z}
\end{equation}
we deduce that $\rho_{\gamma} : G \rightarrow S^{1}$ is a group
morphism. The injectivity of this map results from
Lemma~\ref{lem:RotNumberIdentity}. Since $G$ is compact,
$\rho_{\gamma}$ is also a homeomorphism of $G$ onto a closed
subgroup of $S^{1}$.
\end{proof}

According to Lemma~\ref{lem:RotNumberMonomorphism}, $G$ is either
isomorphic to a finite cyclic group or to $S^{1}$. In the first
case, $f$ is periodic and it is well-known that $f$ is actually
conjugate to a rotation by an angle $2k\pi /n$ about the North-South
axis \cite{CK, Eil,Ker1}. In all of what follows we will suppose
that $G \simeq S^{1}$. In other words, we are given a continuous and
faithful action $\Psi : S^{1} \times S^{2} \rightarrow S^{2}$. We
will establish first:
\begin{lem}
$Fix (G)$ is reduced to two points $N$ and $S$ , $G$ acts freely
on $S^{2} \setminus \set{ N, S }$ and the orbit of every point of
$S^{2} \setminus \set{ N, S }$ is an essential simple closed curve
in $S^{2} \setminus \set{ N, S }$.
\end{lem}

\begin{proof}
Let $\gamma$ be the invariant simple closed curve used to define the
isomorphism $\rho_{\gamma}$. This curve separates the two fixed
points $N$ and $S$ and no element of $G$ other than $Id$ has a fixed
point on $\gamma$. If $G$ has another fixed say $x_{0}$, then the
$G$-orbit of a continuous path joining $x_{0}$ to a point of
$\gamma$ in $S^{2}\ \set{ N ,S }$ gives a homotopy between $\gamma$
and $x_{0}$ in $S^{2} \setminus \set{ N, S }$, contradicting the
fact that $\gamma$ is essential. Hence every $G$-orbit other than
$N$ and $S$ is a simple closed curve and this curve is necessarily
essential in $S^{2} \setminus \set{ N, S }$. Suppose that an element
$g_{0}\in G$ has a fixed point $x_{0}$ in $S^{2} \setminus \set{ N,
S }$ and let $\gamma$($x_{0}$) be the orbit of $x_{0}$ under $G$.
Then $g_{0} = Id$ according to Lemma~\ref{lem:RotNumberIdentity}.
That is, $G$ acts freely on $S^{2} \setminus \set{ N, S }$.
\end{proof}

\begin{lem}\label{lem:MonotoneProperty}
For all $\varepsilon > 0$, there exists $\delta > 0$ so that, if $x$
and $y$ are two distinct points on a $G$-orbit $\gamma$ and $d(x, y)
< \delta$, then at least one of the two arcs delimited by $x$ and
$y$ on $\gamma$ has a diameter less than $\varepsilon$.
\end{lem}

\begin{proof}
Let $\varepsilon > 0$ and $\gamma_{N}$ (resp. $\gamma_{S}$) be a
$G$- orbit in an $\varepsilon$-neighborhood of $N$ (resp. $S$). Let
$A$ be the $G$-invariant annulus bounded by $\gamma_{N}$ and
$\gamma_{S}$. We have only to prove the Lemma for $G$-orbits which
lie in $A$. We can find $\mu > 0$ such that:
\begin{equation}
d(x,g(x)) < \varepsilon /2
\end{equation}
for all $x \in A$ and all $g = \Psi(\theta,.)$ with $\theta \in
I_{\mu} = [-\mu, +\mu]$. Since $G$ acts freely on $A$, there exists
$\delta > 0$ such that:
\begin{equation}
d(x,g(x)) \geq \delta
\end{equation}
for all $x \in A$ and all $g \in S^{1} \ I_{\mu}$. Now, if $x$ and
$y$ are two distinct points on a $G$-orbit $\gamma$ such that $d(x,
y) ) < \delta$ , then $y = \Psi(\theta, x)$ with $\theta \in
I_{\mu}$, and the arc $\Psi ([0, \mu ] , x)$ has diameter less than
$\varepsilon$.
\end{proof}

In order to complete the proof of Theorem~\ref{thm:RegularElliptic},
we are going to show the existence of a ``transversal'' arc to the
$G$-orbit, which will permit us to construct a conjugacy between $G$
and the group of euclidean rotations about the South-North axis (see
also the work of Whitney \cite{Whi} for the existence of transversal
arc-to a family of curves). More precisely:

\begin{cor}\label{cor:TransversalArc}
Given two points $x$ and $y$ which lie on distinct $G$-orbits, there
exists a simple arc $\alpha$ joining $x$ and $y$ which meets each
$G$-orbit in at most one point.
\end{cor}

Let $\gamma$ and $\gamma^{'}$ be two simple closed curves which
separate $N$ and $S$. We write $\gamma \leq \gamma^{'}$ (resp.
$\gamma < \gamma^{'}$) iff $\gamma$ is contained in the closed
(resp. open) disc bounded by $\gamma^{'}$ and containing $S$. This
relation induces a total order on the set of $G$-orbits (with the
convention that $S \leq \gamma$ and $\gamma \leq N$ for all
$G$-orbit $\gamma$).

For any point $x \in S^{2}$, we let $\gamma(x)$ be the $G$-orbit
of $x$. With this notation, we have the following definition. A
finite collection of indexed points $\set{x_{0}, x_{l},\dotsc,
x_{n} }$ such that $d(x_{k}, x_{k+l}) < \mu$ and $\gamma(x_{k}) <
\gamma (x_{k+l})$ is called a \textit{monotone $\mu$-chain} from
$x_{0}$ to $x_{n}$. We will establish first:

\begin{lem}\label{lem:MonotoneChain}
For all $\varepsilon > 0$, there exists $\delta > 0$ so that two
points $x$ and $y$ with $d(x, y) < \delta$ and $\gamma(x) <
\gamma(y)$ can be joined by a monotone $\mu$-chain of diameter less
than $\varepsilon$ for all $\mu > 0$.
\end{lem}

\begin{proof}
Let $\varepsilon > 0$ and $\delta > 0$ as in
Lemma~\ref{lem:MonotoneProperty}. Suppose that $x,y$ are two points
lying on distinct $G$-orbits and such that $d(x, y) < \delta/2$.
Given $\mu > 0$ ($\mu < \delta$), we can find a finite sequence of
$G$-orbits
\begin{equation}
\gamma_{0} = \gamma (x),\gamma_{1},\dotsc ,\gamma_{n} =\gamma(y)
\end{equation}
such that $d_{H}(\gamma_{k}, \gamma_{k+l}) < \mu/3$.

Let $\overline{xy}$ be the geodesic arc connecting $x$ and $y$ whose
length is less than $\delta/2$. This arc meets each intermediate
curve $\gamma_{k}$. Hence, we can choose, for each $k$, a point
$x_{k}$ of $\overline{xy}$ on $\gamma_{k}$. If $d(x_{k}, x_{k+l}) <
\mu$ for all $k$, we are done. If not, we are going to show that for
each pair $\set{x_{r},x_{r+1}, }$ for which $d(x_{r}, x_{r+1}) \geq
\mu$, we can join $x_{r}$ and $x_{r+1}$ by a monotone $\mu$-chain
which lies in a $2\varepsilon$-neighborhood of $\overline{xy}$. By
joining together such chains, we obtain a monotone $\mu$-chain from
$x$ to $y$ of diameter less than $4\varepsilon + \delta/2$.

Choose a point $x^{'}_{r+1}$ on $\gamma_{r+1}$ such that $d(x_{r},
x^{'}_{r+1}) < \mu /3$. Hence, $d(x_{r+1}, x^{'}_{r+1}) < \delta$
and one of the arcs delimited by these two points on $\gamma_{r+1}$,
$x_{r+1}x^{'}_{r+1}$, has a diameter less than $\varepsilon$.

Divide this arc into $s$ subarcs of diameter less than $\mu/3$ and
call the intermediate points
\begin{equation}
x^{'}_{r+1} = z^{0}_{r+1},z^{1}_{r+1}\ldots ,z^{s}_{r+1}= x_{r+1}
\end{equation}
We choose some curves
\begin{equation}
    \gamma_{r} =  \gamma_{0} <  \gamma_{1} < \dotsb <  \gamma_{s} =  \gamma_{r+1}
\end{equation}
For each $k \in \set{1, \ldots, s - 1  }$, $dH(\gamma^{k},
\gamma_{r+l}) < \mu/3$ and we can choose on each intermediate
curve $\gamma^{k}$ a point $x^{k}_{r+1}$ such that $d(x^{k}_{r+1},
z^{k}_{r+1}) < \mu/3$. Hence
\begin{equation}
x^{0}_{r+1} = x_{r}, x^{1}_{r+1},\ldots, z^{s}_{r+1} = x_{r+1}
\end{equation}
is a monotone $\mu$-chain.
\end{proof}

\begin{proof}[Proof of Lemma 4.12.]
From Lemma~\ref{lem:MonotoneChain}, we can choose a sequence of
numbers $\delta_{n} > 0$ so that any two points with distance $d(x,
y) < \delta_{n}$, can be joined, for all $\mu > 0$ by a monotone
$\mu$-chain of diameter less than $1/2^{n}$.

We start by choosing a monotone $\delta_{0}$-chain $X_{0}$ from $x$
to $y$. Inductively, once $X_{n}$ has been defined, we join each
consecutive pair ${x^{n}_{k}, x^{n}_{k+1}}$ of $X_{n}$, by a
monotone $\delta_{n+1}$-chain of diameter less than $1/2^{n+1}$ to
obtain $X_{n+1}$ and we set
\begin{equation}
X = \bigcup_{n\in\mathbb{N}}X_{n}.
\end{equation}
It is then a standard fact (\cite[Theorem 2.27]{HY}) that the
closure $\overline{X}$ of $X$ in $S^{2}$ is a simple arc joining $x$
and $y$ with the required property.
\end{proof}

To complete the proof of Theorem~\ref{thm:RegularElliptic}, we
choose an arc $\alpha$ given by Corollary~\ref{cor:TransversalArc}
from $N$ to $S$ , and we let $x(r), r \in [0,+\infty]$ be a
parametrization of this arc. The map
\begin{equation}
h : re^{i\theta} \mapsto \Psi(\theta, x (r) )
\end{equation}
is the required conjugacy between $G$ and the group of euclidean
rotations about the $z$-axis.


\subsection{Orientation-reversing homeomorphisms of the
sphere}\label{subsec:SphereNegative}

Every orientation-reversing conformal automorphism of the Riemann
Sphere $\widehat{\mathbb{C}}$ is of the form
\begin{equation}
    f(z) = (a\overline{z} + b)/(c\overline{z} + d),
\end{equation}
where the coefficients satisfy $ad - bc \neq 0 $. We call them fractional
reflections following Maskit \cite{Mas}. The fixed point set of a
fractional reflection is either empty, one point, two points or a
circle in $\widehat{\mathbb{C}}$. Fractional reflections are classified by the number
of fixed points:
\begin{itemize}
    \item A transformation $f$ with a circle of fixed points is a reflection. It is topologically conjugate to the map $z\mapsto \overline{z}$.
    \item A transformation with exactly two fixed points is
semi-hyperbolic, It is topologically conjugate to the map $z\mapsto 2\overline{z}$.
    \item A transformation with exactly one fixed point is
semi-parabolic. It is topologically conjugate to the map $z\mapsto \overline{z} +1$.
    \item A transformation with no fixed points is semi-elliptic. It
is topologically conjugate to the map $z\mapsto e^{i\theta}/ \overline{z}$, with $\theta \neq 0$ mod $2\pi$.
   \end{itemize}

\begin{thm}\label{thm:SphereNegative}
An orientation-reversing homeomorphism of the sphere $S^{2}$ with a
totally disconnected singular set is topologically conjugate to a
fractional reflection.
\end{thm}

\begin{proof}
1) Suppose first that the singular $\Sigma (f)$ is empty. Then
$\Sigma (f^{2}) = \Sigma (f) = \emptyset$ and according to
Theorem~\ref{thm:MainSphere}, we may suppose after a topological
conjugacy if necessary that $f^{2}$ is a euclidean rotation around
the vertical axis.

If $f$ has a fixed point, there are arbitrarily small closed discs
invariant under $f$ (c.$f$.
Corollary~\ref{cor:FixedPointNeighborhood}). On the boundary of each
of these discs, $f$ reverses the order and has therefore two fixed
points. In particular, $f$ and hence $f^{2}$ have an infinite number
of fixed points and therefore $f^{2} = id$. In that case $f$ is
conjugate to a reflection \cite{CK}.

If $f$ has no fixed points, $f$ permutes the two fixed points of
$f^{2}.$ If $f^{2}$ is a periodic rotation, we refer also to
\cite{CK} for a proof that $f$ is conjugate to the map $z \mapsto
e^{2i\pi p/q}/\overline{z}$. If $f^{2}$ is an irrational rotation by
angle $2\alpha$, the family of simple closed curves invariant under
$f^{2}$ is unique. Hence, $f$ permutes these curves and induces a
continuous, reversing-order involution on $S^{2}/G \cong [0,1]$
where $G$ is the closure of the group generated by $f^{2}$.
Therefore one and exactly one of theses curves is invariant under
$f$. This curve divides the sphere into exactly two discs which are
permuted by $f$ and each of them contains a fixed point of $f^{2}$.
In one of these discs, we choose a ``transverse'' arc which joins
the fixed point of $f^{2}$ to the boundary of this disc. We map this
arc onto an arc with similar properties relatively to the map $z
\mapsto e^{2i\pi\alpha}/\overline{z}$ and extend this map by
iteration under $f$. This map is clearly well-defined and extends
into a topological conjugacy between $f$ and the map $z \mapsto
e^{2i\pi\alpha}/\overline{z}$.

2) If $\Sigma (f)$ contains exactly one point $N$ then $f(N) = N$
and the group generated by $f$ acts freely and properly on
$S^{2}\setminus \set{ N  }$. The quotient space $\Gamma_{f}$ is
homeomorphic to the open M\"{o}bius strip and we can show as in section
4.1, that $f$ is conjugate to the map $z  \mapsto \overline{z} + 1$.

3) If $\Sigma (f)$ contains exactly two points $N$ and $S$, then
$f(N) = N$  and $f(S)= S$. The group generated by $f$ acts freely
and properly on the cylinder $I = S^{2} \setminus \set{ N ,S }$. The
quotient space $\Gamma_{f}$ is homeomorphic to the Klein bottle and
we can show as in section~\ref{subsec:Hyperbolic}, that $f$ is
conjugate to the map $z \mapsto 2\overline{z}$.
\end{proof}

We conclude this paragraph by stating the analogous result for the
closed disc $\mathbb{D}^{2}$. Let $f$ be a homeomorphism of
$\mathbb{D}^{2}$ whose singular set $\Sigma (f)$ is totally
disconnected. The same property holds for the homeomorphism induced
on the sphere $S^{2}$ viewed as the double of the disc. According to
Theorem~\ref{thm:MainSphere}, $f$ has at most two singular points
which lie necessarily on the boundary of $\mathbb{D}^{2}$. The
proofs which have been given for the sphere may be adapted in the
case of the disc without introducing any new subtlety to establish
the following:
\begin{thm}
A homeomorphism of the closed disc with a totally disconnected
singular set is topologically conjugate to a fractional linear
transformation or a fractional reflection of the disc according to
whether it is orientation-preserving or orientation-reversing.
\end{thm}


\subsection{The projective plane}\label{subsec:Projective}

Let $f$ be a homeomorphism of the projective plane
$\mathbb{P}^{2}(\mathbb{R})$. If $\Sigma(f)$ is totally
discontinuous, $\Sigma(f)$ must be empty according to
Corollary~\ref{cor:SingularSetEmpty}. Let $\widetilde{f}$ be the
unique orientation-preserving lift of $f$ to the universal cover
$S^{2}$ of $\mathbb{P}^{2}(\mathbb{R})$. Then $\widetilde{f}$ is a
regular homeomorphism which commutes with the covering involution
$s:z \mapsto -1/\overline{z}$.

If $\widetilde{f}$ is non periodic, the closure of the group
generated by $\widetilde{f}$, $G$, is isomorphic to $S^{1}$ and
there is a unique family of simple closed curves invariant under
$\widetilde{f}$. Since $s$ commutes with $\widetilde{f}$, $s$
permutes these curves and induces a continuous,
orientation-reversing involution of $S^{2}/G \cong [0, 1]$. One and
exactly one of theses curves is invariant under $s$. We leave it to
the reader to show that in that case, the conjugacy between
$\widetilde{f}$ and the standard rotation $z\mapsto e^{2i \pi \theta
z}$ may be chosen to commute with $s$. In other words, $f$ and the
standard map induced on $\mathbb{P}^{2}(\mathbb{R})$ by $z \mapsto
e^{2i \pi \theta z}$ are conjugate on $\mathbb{P}^{2}(\mathbb{R})$.

If $\widetilde{f}$ is periodic, then $\widetilde{f}$ is conjugate to
a periodic rotation. The quotient space $S^{2}/ \widetilde{f}$ is a
sphere and $s$ induces an orientation-reversing involution of this
sphere. Therefore, we can find a simple closed curve on $S^{2}$,
which separates the two fixed points of $\widetilde{f}$ and which is
invariant both by $\widetilde{f}$ and $s$. These considerations may
be used to show that in that case also an equivariant conjugacy
between $\widetilde{f}$ and $z \mapsto e^{2i \pi p/q z}$ can be
constructed. We have finally proven the following:

\begin{thm}
Let $f$ be a regular homeomorphism of the projective plane
$\mathbb{P}^{2}(\mathbb{R})$. Then $f$ is topologically conjugate to
a standard rotation of the projective plane.
\end{thm}


\section{The case of surfaces with $\chi(M^{2}) \leq
0$}\label{sec:NegativeEuler}

Let $f$ be a homeomorphism of a compact surface $M^{2}$ with
$\chi(M^{2}) \leq 0$. If $\Sigma (f)$ contains no non-degenerate
continuum, $f$ is regular everywhere according to
Corollary~\ref{cor:SingularSetEmpty}.


\subsection{General results on regular homeomorphisms of
surfaces}\label{subsec:GeneralRegular}

Let $M^{2}$ be a closed orientable surface of genus $g \geq 1$ and
$\pi : \widetilde{M}^{2} \rightarrow M^{2}$ the universal cover of
$M^{2}$. We can identify $\widetilde{M}^{2}$ either to the euclidean
plane $\mathbb{R}^{2}$ or to the Poincar\'{e} disc $\mathbb{D}$ in such
a way that $M^{2}$ is homeomorphic to the quotient of
$\widetilde{M}^{2}$ by a discrete subgroup $\Gamma$ of euclidean
translations or hyperbolic isometries according to whether
$\widetilde{M}^{2}$ is $\mathbb{R}^{2}$ or $\mathbb{D}$. The metric
we shall use on $M^{2}$ is the quotient metric defined on
$\widetilde{M}^{2}/\Gamma$ by:
\begin{equation}
d(\pi(\widetilde{x}), \pi(\widetilde{y})) = \inf_{g,h \in \Gamma} \widetilde{d}(g.\widetilde{x}, h.\widetilde{y})
\end{equation}
where $\widetilde{d}$ is the natural metric on $\widetilde{M}^{2}$.

There is another metric on $\widetilde{M}^{2}$ that we shall use in
the following, namely the spherical metric. The Alexandroff
compactification of $\widetilde{M}^{2}$, $\widetilde{M}^{2}\cup
{\infty}$, is homeomorphic to the sphere $S^{2}$. Hence, the
standard metric of $S^{2}$ induces a metric that we shall call
$\partial$ on $\widetilde{M}^{2}$. These two metrics $\widetilde{d}$
and $\partial$ are not uniformly equivalent on $\widetilde{M}^{2}$
but $Id : (\widetilde{M}^{2}, \widetilde{d}) \rightarrow
(\widetilde{M}^{2},\partial)$ is uniformly continuous.

\begin{lem}\label{lem:RegularLift}
Let $f$ be a regular homeomorphism of a closed orientable surface
$M^{2}$ of genus $g \geq 1$ and $\widetilde{f}$ be any lift of $f$
to the universal cover $\widetilde{M}^{2}$ of $M^{2}$. If
$\overline{f}$ is the unique continuous extension of $\widetilde{f}$
to $S^{2} = \widetilde{M}^{2} \cup \set{ \infty
 }$ then $\Sigma (\overline{f})\subset \set{ \infty  }$.
\end{lem}

\begin{proof}
Let $\widetilde{f}$ be any lift of $f$ on $\widetilde{M}^{2}$.
Every point $\widetilde{x}$ is regular under $\widetilde{f}$ for
the metric $\widetilde{d}$. Let $\partial$ be the standard metric
on $S^{2}$. Then, $Id : (\widetilde{M}^{2}, \partial) \rightarrow
(\widetilde{M}^{2}, \widetilde{d})$ is continuous and $Id :
(\widetilde{M}^{2}, \widetilde{d}) \rightarrow (\widetilde{M}^{2},
\partial)$ is uniformly continuous. Therefore, a point
$\widetilde{x}$ which is regular for $(\widetilde{f} ,
\widetilde{d})$ is regular for $(\overline{f},\partial )$ and
hence $\Sigma (\overline{f})\subset \set{ \infty  }$.
\end{proof}

\begin{cor}\label{cor:RegularPositiveGenius}
A regular homeomorphism of a closed orientable surface of genus $g
\geq 1$ which is homotopic to the identity and which has a fixed
point is the identity.
\end{cor}

\begin{proof}
Let $\widetilde{f}$ be any lift of $f$ on $\widetilde{M}^{2}$ which
has a fixed point, Since $f$ is homotopic to the identity,
$\widetilde{f}$ commutes with all covering translations and has
therefore an infinite number of fixed points. According to
Lemma~\ref{lem:RegularLift} and Theorem~\ref{thm:MainSphere},
$\overline{f}$ and hence $\widetilde{f}$ must be equal to the
identity, which completes the proof.
\end{proof}

\begin{thm}\label{thm:MainNegativeEuler}
A regular homeomorphism of a compact surface of negative Euler
characteristic is periodic.
\end{thm}

\begin{proof}
If the boundary of $M^{2}$ is not empty, the natural extension of
$f$ to the double $DM^{2}$ of $M^{2}$ is still regular and since
$\chi(DM^{2}) = 2\chi(M^{2})$, we are reduced to prove
Theorem~\ref{thm:MainNegativeEuler} for closed surfaces. Moreover,
by passing to the orientation covering of $M^{2}$ and by considering
$f^{2}$ instead of $f$ if necessary, we may assume that $M^{2}$ is
orientable and that $f$ is orientation-preserving. So, let $f$ be a
regular orientation-preserving homeomorphism of a closed orientable
surface. Since $f$ is recurrent we can find a positive integer n
such that $f^{n}$ is homotopic to the identity and since
$\chi(M^{2}) < 0$, Lefschetz's formula implies that $f^{n}$ has a
fixed point. According to Corollary~\ref{cor:RegularPositiveGenius},
$f^{n}$ is equal to the identity, which completes the proof.
\end{proof}


\subsection{Orientation-preserving homeomorphisms of the
torus}\label{subsec:TorusPositive}

The \textit{translations} of the torus are the maps of the torus
induced by standard translations $\tau_{\alpha,\beta} : (s, t)
\mapsto (s + \alpha, t + \beta)$ of the plane. Each of these maps
is a regular transformation of the torus and it is periodic if and
only if $(\alpha,\beta) \in \mathbb{Q}^{2}$. Moreover, two such
maps $\tau_{\alpha,\beta}$ and $\tau_{\gamma,\delta}$ are
topologically conjugate if and only if the two vectors $(\alpha,
\beta)$ and $(\gamma, \delta)$ can be mapped one onto the other by
a matrix $A \in GL(2, \mathbb{Z})$. The aim of this paragraph is
to establish the following:

\begin{thm}\label{thm:MainTorusPositive}
A regular, orientation-preserving and non periodic homeomorphism of
the torus $\mathbb{T}^{2}$ is topologically conjugate to a non
periodic translation of the torus.
\end{thm}

\begin{rem}
A complete classification of periodic transformations of the torus
has been given by Brouwer \cite{Bro2, Bro3, Yok}. We will not give
here the list of them which can be find in these references.
\end{rem}

Let $f$ be a regular orientation-preserving homeomorphism of the
torus and let $A \in SL(2, \mathbb{Z})$ be the induced matrix on
$\pi_{1}(\mathbb{T}^{2}) \cong H_{l}(\mathbb{T}^{2})\cong
\mathbb{Z}^{2}$. Since $f$ is recurrent, we get that $A^{n} = Id$
for some $n > 0$.

If $A \neq Id$ then the Lefschetz number of $f$, $L(f) = 2 - Tr(A)
\neq 0$ and $f$ has a fixed point. Hence, $f^{n} = Id$ according to
Corollary~\ref{cor:RegularPositiveGenius}.

If $A = Id$, then for any lift $\widetilde{f}$ of on
$\mathbb{R}^{2}$ the following relation holds:

\begin{equation}
\widetilde{f}(\widetilde{x} + v) = \widetilde{f} (\widetilde{x}) +
v,\qquad \forall \widetilde{x} \in \mathbb{R}^{2}, \forall v \in
\mathbb{Z}^{2}
\end{equation}
in other words $f$ commutes with integer translations and
$\widetilde{f} (\widetilde{x}) - \widetilde{x}$ is uniformly
bounded on $\mathbb{R}^{2}$.

\begin{lem}\label{lem:RotVector}
Let $f$ be a regular and orientation-preserving homeomorphism of the
torus which acts trivially on $\pi_{1}(\mathbb{T}^{2})$ and let $f$
be any lift of $f$ on $\mathbb{R}^{2}$. Then
\begin{enumerate}
    \item $\theta(\widetilde{f}, \widetilde{x}) = \lim((\widetilde{f}^{n}(\widetilde{x}) - \widetilde{x})/n)$ exists and is independent of $\widetilde{x}$. We shall call it the translation vector of $\widetilde{f}$ and denote it by $\theta (\widetilde{f})$.
    \item $\widetilde{f} = (0, 0)$ iff $\widetilde{f}$ has a fixed point.
\end{enumerate}
\end{lem}

\begin{proof}
1) Let $K_{n}$ be the closure of the set
\begin{equation}
    \set{ \frac{\widetilde{f}^{m}(\widetilde{x})-\widetilde{x}}{m};\quad m\geq n,\widetilde{x}\in\mathbb{R}^{2} }.
\end{equation}
Since the map $\widetilde{f}(\widetilde{x}) - \widetilde{x}$ is
bounded on $\mathbb{R}^{2}$, $(K_{n})$ is a decreasing sequence of
compact sets. We are going to show that $diam \,(K_{n})
\rightarrow 0$ as $n \rightarrow \infty$. Let $\varepsilon > 0$.
Since $f$ is regular and thus recurrent, we can find a positive
integer $r$ and an integer vector $v$ such that:

\begin{equation}
    \norm{ \widetilde{f}^{r}(\widetilde{x})-(\widetilde{x}+v) } < \varepsilon,\quad \forall\widetilde{x}\in\mathbb{R}^{2}
\end{equation}
and therefore
\begin{equation}
    \norm{\frac{\widetilde{f}^{kr}(\widetilde{x)}-\widetilde{x}}{kr}-\frac{v}{r} } < \frac{\varepsilon}{r},\quad \forall k >0, \forall\widetilde{x} \in \mathbb{R}^{2}.
\end{equation}
Let $M$ be a bound of $\norm{ \widetilde{f}(\widetilde{x}) -
\widetilde{x} }$ on $\mathbb{R}^{2}$. We have then:
\begin{equation}
    \norm{\frac{\widetilde{f}^{n}(\widetilde{x)}-\widetilde{x}}{n}-\frac{v}{r} } < \frac{\varepsilon}{r}+\frac{2rM}{n},\quad \forall n \geq r, \forall\widetilde{x} \in \mathbb{R}^{2}
\end{equation}
which shows that $diam \,(K_{n}) < 2\varepsilon$ for $n$ big
enough. Therefore: $\widetilde{f}^{n} (\widetilde{x}) -
\widetilde{x}$ converges uniformly on $\mathbb{R}^{2}$ to some
constant vector $\theta(\widetilde{f})$ as $n \rightarrow
+\infty$.

2) If $f$ has a fixed point, then clearly $\theta(\widetilde{f})=0$.
Conversely, if $f$ has no fixed point, then according to
Lemma~\ref{lem:RegularLift} and Theorem~\ref{thm:MainSphere}, $f$
has only wandering points.  Hence there exists $n_{0}$ such that for
every $n \geq n_{0}$, $\norm{ \widetilde{f}^{n}(0, 0)  } > 1$.
Choose a positive number $\varepsilon < 1/2$, and let $n_{1} >
n_{0}$ and $v \in \mathbb{R}^{2}$ such that:
\begin{equation}
    \norm{\widetilde{f}^{n_{1}}(\widetilde{x})-(\widetilde{x}+v) }<\varepsilon,\quad \forall\widetilde{x}\in\mathbb{R}^{2}.
\end{equation}

This inequality shows first that $v$ cannot be the vector $(0,
0)$, since in, $\widetilde{f}^{n_{1}}(0, 0)$ is not contained in
the ball of center $(0, 0)$ and of radius $\varepsilon$. Then, we
have:
\begin{equation}
    \norm{\frac{\widetilde{f}^{kn_{1}}(\widetilde{x)}-\widetilde{x}}{kn_{1}}-\frac{v}{n_{1}} } < \frac{\varepsilon}{n_{1}},\quad \forall k >0, \forall\widetilde{x} \in \mathbb{R}^{2}
\end{equation}
and letting $k\rightarrow +\infty$, we obtain:
\begin{equation}
    \norm{\theta (\widetilde{f})-\frac{v}{n_{1}} } \leq \frac{\varepsilon}{n_{1}},
\end{equation}
which shows that $\theta (\widetilde{f})$ cannot be zero.
\end{proof}

\begin{rem}
This Lemma is still true if we replace the statement $f$ regular by
$f$ recurrent which is a weaker hypothesis. In fact, the proof of
the first part of the Lemma does indeed only use this hypothesis. We
can also remark that 1) says precisely that the rotation set defined
by Misiurewicz and Zieman in \cite{MZ} is reduced to a point for a
recurrent homeomorphism of the torus. To prove 2), we need to know
the fact that an orientation-preserving and fixed point free
homeomorphism of the plane has only wandering points which is a
corollary of Brouwer's Lemma on translation arcs \cite{Fat, Gui}
which we have not used here. 2) can be considered as a particular
case of a result of Franks \cite{Fra}.
\end{rem}

Let $\widetilde{g}$ be another lift of $f$. Then
$\theta(\widetilde{f}) - \theta(\widetilde{g}) \in \mathbb{Z}^{2}$
and  the class of $\theta(\widetilde{f})$ modulo $\mathbb{Z}^{2}$
is independent of the particular choice of the lift
$\widetilde{f}$. We shall call it the \textit{rotation vector} of
$f$ and denote it by $\rho (f)$.

As in Theorem~\ref{thm:Equicontinuity}, let $G$ be the closure of
the family $\set{ f^{n} ; \, n \in \mathbb{Z} }$. $G$ is a
commutative compact group of regular homeomorphisms of
$\widetilde{T}^{2}$ and each element of $G$ acts trivially on the
fundamental group of $\widetilde{T}^{2}$. Moreover, we have:

\begin{lem}
Let $f$ be a regular homeomorphism of the torus
$\widetilde{T}^{2}$ acting trivially on the fundamental group, and
let $G$ be the closure of the group generated by $f$. Then the map
$p : G \rightarrow \mathbb{T}^{2}, g \mapsto \rho(g)$ induces a
bicontinuous isomorphism from $G$ onto a compact subgroup of
$\mathbb{T}^{2}$.
\end{lem}

\begin{proof}
We will show first that $p$ is a group morphism. Let $g, h \in G$
and let $\widetilde{g}, \widetilde{h}$ be lifts of $g$ and $h$
respectively. Since $g$ and $h$ commute, there exists $v \in
\mathbb{Z}^{2}$ such that
\begin{equation}
(\widetilde{h}^{-1} \circ \widetilde{g} \circ \widetilde{h})- \widetilde{g} = v.
\end{equation}
Since $\widetilde{h}$ and $\widetilde{g}$ commute with the integer
translations and since other lifts are obtained by composing
$\widetilde{g}$ and $\widetilde{h}$ with integer translations, one
shows easily that the vector $v$ does not depend on the lifts
$\widetilde{g}$ and $\widetilde{h}$. Furthermore, from the
commutativity of $\widetilde{g}$ and $\widetilde{h}$ with the
integer translations, one obtains inductively
\begin{equation}
\widetilde{h}^{-n} \circ \widetilde{g} \circ \widetilde{h}^{n} - \widetilde{g} = nv
\end{equation}
and this relation also holds for any other lift of $g$ and
$h^{n}$. Since $h$ is recurrent, there exists $n > 0$ such that
$h^{n}$ is arbitrarily close to identity. Hence, we can find a
lift of $h^{n}$ close to identity and therefore $v = 0$; that is
$\widetilde{g}$ and $\widetilde{h}$ commute. Hence, we can write:
\begin{equation}
\frac{(\widetilde{g}\circ\widetilde{h})^{n}(\widetilde{x})-\widetilde{x}}{n}
= \frac{\widetilde{g}^{n}
(\widetilde{h}^{n}(\widetilde{x}))-\widetilde{h}^{n}(\widetilde{x})}{n}
+ \frac{\widetilde{h}^{n}(\widetilde{x})-\widetilde{x}}{n}
\end{equation}
and since $(\widetilde{g}^{n}(x) - x)/n$ converges uniformly to
$\theta (\widetilde{g})$ on $\mathbb{R}^{2}$, we get
\begin{equation} \theta(\widetilde{g} \circ \widetilde{h}) =
\theta(\widetilde{g}) + \theta(\widetilde{h}).
\end{equation}
That is $\rho: G \rightarrow \mathbb{T}^{2}$ is a group morphism.
According to Lemma~\ref{lem:RotVector} and
Corollary~\ref{cor:RegularPositiveGenius}, this morphism is
necessarily injective. The continuity of $\rho$ results from the
fact that
\begin{equation}
    \norm{ \widetilde{g}-Id  } <\varepsilon \Rightarrow \norm{ \theta (\widetilde{g}) } \leq \varepsilon.
\end{equation}
\end{proof}

A compact subgroup of $\mathbb{T}^{2}$ is either a finite group,
$\mathbb{T}^{2}$ or the product of $S^{1}$ by a finite cyclic
group \cite{Bou}. In the first case $f$ is periodic. We shall
leave aside the third case until the end of this section.

\begin{lem}
Assume that $\rho(G) = \mathbb{T}^{2}$. Then there is a
topological conjugacy between $G$ and the group of translations of
$\mathbb{T}^{2}$.
\end{lem}

\begin{proof}
Fix a point $x_{0}$ in $\mathbb{T}^{2}$ and let $\phi:
\mathbb{T}^{2} \rightarrow \mathbb{T}^{2}$ be defined as follows.
For each $t \in \mathbb{T}^{2}$ there is a unique $g_{t} \in G$
with $\rho(g_{t}) = t$ and we define a continuous map $\phi$ by
the following:
\begin{equation}
\phi(t) = g_{t}(x_{0}).
\end{equation}
If $\phi(t) = \phi(s)$, then $(g_{s}^{-1} \circ g_{t}) (x_{0}) =
x_{0}$ and so $g_{s}^{-1} \circ g_{t}$ is the identity map
(Corollary~\ref{cor:RegularPositiveGenius}). As $\rho$ is an
isomorphism this implies that $s = t$. Hence $\phi$ is one-to-one
and is thus a homeomorphism of $\mathbb{T}^{2}$ onto
$\phi(\mathbb{T}^{2}) = \mathbb{T}^{2}$, by the invariance of
domain.

It remains to be shown that $\phi$ is a conjugacy between the
group of translations and $G$. Given $s \in \mathbb{T}^{2}$, we
denote by $\tau_{s}$ the translation by $s$. For every $t \in
\mathbb{T}^{2}$, we have
\begin{equation}
\phi \circ \tau_{s}(t) = \phi (t + s) = g_{s+t}(x_{0}) = g_{s}\circ g_{t}(x_{0}) =g_{s}(\phi(t))
\end{equation}
which completes the proof.
\end{proof}

\begin{rem}
We have proved in fact that every faithful and continuous action
of $\mathbb{T}^{2}$ on $\mathbb{T}^{2}$ is isomorphic to the
standard action.
\end{rem}

From now on, we shall assume that $\rho(G)$ is (up to a linear conjugacy of $\mathbb{T}^{2}$) the subgroup
\begin{equation}
\left ( \frac{1}{q}\cdot \mathbb{Z} \right)/\mathbb{Z}\times S^{1}\subset \mathbb{T}^{2}.
\end{equation}
Let g be the unique element of $G$ whose rotation vector is $(1/q,
0)$ and let $G_{0}$ be the connected component of the identity. $G$
is clearly the direct product of $G_{0}$ by $< g >$, the finite
group (isomorphic to $\mathbb{Z}/q\mathbb{Z}$) generated by $g$.
Each $G$-orbit is a family of $q$ distinct simple closed curves
since no element of $G$ has a fixed point. These curves divide the
torus into $q$ distinct topological annulus $A_{0}, A_{1}, ...,
A_{q-1}$ which are permuted by $g$ and each one of these annulus is
invariant under $G_{0}$. According to the results of
section~\ref{sec:PositiveEuler}, the restriction of the action of
$G_{0}$ on each annuli $A_{i}$ is conjugate to the standard action
of $S^{1}$. From these considerations, we deduce the following
results which complete the proof of
Theorem~\ref{thm:MainTorusPositive}.

\begin{lem}
There exists a simple path $\sigma_{0}: [0, 1/q] \rightarrow
\mathbb{T}^{2}$ joining $(0, 0)$ to $g ((0, 0))$ such that for any
$s, t \in [0, 1/q], s \neq t$ and $\set{s, t  } \neq \set{0, 1/q
}$, $\sigma_{0}(s)$ and $\sigma_{0}(t)$ are on distinct
$G_{0}$-orbits.
\end{lem}

\begin{cor}\label{cor:TorusDetails}
We have the following:
\begin{enumerate}
    \item Let $\sigma: [0, 1] \rightarrow \mathbb{T}^{2}$ be defined by $\sigma(\frac{i}{q} + t) = g^{i}(\sigma_{0}(t)),t \in [0, 1/q]$. Then $\sigma$ is a simple closed curve invariant under $g$, which meets each orbit of $G_{0}$ in exactly one point.
    \item Let $\phi: \mathbb{T}^{2} \rightarrow \mathbb{T}^{2}$ be the map defined by $\phi(s, t ) ) = g_{(0,t)} (\sigma (s))$ where $_{(0,t)}$ is the unique element of $G$ whose rotation vector is equal to $(0, t)$. Then $\phi$ is a conjugacy between $G$ and the subgroup of translations by elements of $( \frac{1}{q}\cdot \mathbb{Z} )/\mathbb{Z}\times S^{1}\subset \mathbb{T}^{2}$.
   \end{enumerate}
\end{cor}


\subsection{Orientation-reversing homeomorphisms of the
torus}\label{subsec:TorusNegative}

Let $f$ be an orien\-tation-reversing homeomorphism of the torus
$\mathbb{T}^{2}$. According to Theorem~\ref{thm:mainthm}(2), $f^{2}$
is periodic or is conjugate to a translation by a vector with at
least one irrational coordinate. The two families of maps $(s,t )
\mapsto (-s, t + \alpha)$ and $(s,t ) \mapsto (-s, t + s + \alpha)$,
where $\alpha \in (\mathbb{R} \setminus \mathbb{Q})/\mathbb{Z}$ are
standard examples of regular, orientation-reversing and non periodic
homeomorphism of the torus $\mathbb{T}^{2}$. Two such maps are not
conjugate. We shall show in this paragraph that, up to topological
conjugacy, there are no other such map, More precisely:

\begin{thm} Let $f$ be a regular,
orientation-reversing and non periodic homeomorphism of the torus
$\mathbb{T}^{2}$. Then, there exists $\alpha \in ( \mathbb{R} \ \mathbb{Q})/\mathbb{Z}$ such that $f$ is topologically
conjugate either to $(s,t) \mapsto (-s, t + \alpha)$ or $(s, t ) \mapsto (-s, t
+ s + \alpha)$.
\end{thm}

\begin{proof}
Let $A \in GL(2, \mathbb{Z})$ be the induced matrix on
$\pi_{1}(\mathbb{T}^{2})\cong H_{l}(\mathbb{T}^{2})\cong
\mathbb{Z}^{2}$. For any lift $\widetilde{f}$ of $f$ to the
universal cover $\mathbb{R}^{2}$ of $\mathbb{T}^{2}$ we have:
\begin{equation}
\widetilde{f}(\widetilde{x} + v) = \widetilde{f}(\widetilde{x}) + A \cdot v,
\end{equation}
for every integer vector $v$. Furthermore, since $f$ is not
periodic, $A^{2} = I$ according to the results of
section~\ref{subsec:TorusPositive}. Let $\theta(\mathbb{f}^{2},
\widetilde{x})$ be the function defined in
Lemma~\ref{lem:RotVector}. An easy computation shows that:
\begin{equation}
    \theta (\widetilde{f} \circ \widetilde{f}^{2} \circ \widetilde{f}^{-1},\widetilde{f}(\widetilde{x}) = A \cdot \theta (\widetilde{f}^{2},\widetilde{x})
\end{equation}
since the function $\widetilde{f}(\widetilde{x}) -
A(\widetilde{x})$ is bounded on $\mathbb{R}^{2}$. Moreover, since
$\theta(\widetilde{f}^{2}, \widetilde{x})$ has been shown to be
independent of $\widetilde{x}$, we get:
\begin{equation}
A\cdot\theta(\widetilde{f}^{2}) =\theta(\widetilde{f}^{2}).
\end{equation}
In other words $\theta(\widetilde{f}^{2})$ is an eigenvector of
$A$. Since $f$ reverses orientation, $\det A = -1$ and the
eigenvalues of $A$ are $+1$ and $-1$. In $GL(2, \mathbb{Z})$ there
are two conjugacy classes of such matrices
\begin{equation}
    \begin{array}{cc}
      \left(%
            \begin{array}{cc}
                 -1 & 0 \\
                  0 & 1 \\
            \end{array}%
      \right) & (\text{Type 1})\\
     \end{array}
\end{equation}
and
\begin{equation}
    \begin{array}{cc}
      \left(%
            \begin{array}{cc}
                 -1 & 0 \\
                  1 & 1 \\
            \end{array}%
      \right) & (\text{Type 2})\\
     \end{array}
\end{equation}
In both cases, we can find a new system of coordinates on the
torus in which we have $\theta(\widetilde{f}^{2}) = (0, 2\alpha)$.

Since $f^{2}$ is not periodic
by assumption, $\alpha$ must be an irrational number and the closure, $G$,
of the group generated by $f^{2}$ is isomorphic to $S^{1}$. Each $G$-orbit is
a simple closed curve whose homology class is $(0, 1)$ in this
system of coordinates.

According to Corollary~\ref{cor:TorusDetails}, we can suppose, after
a change of coordinates if necessary, that $f^{2}$ is the map
defined by $(s,t ) \mapsto (s, t + 2\alpha)$. The circles $\set{ s }
\times S^{1}$ become then the closures of $f^{2}$-orbits. In
particular, $f$ maps each such circle onto another one and its
analytic expression must be $(s, t)\mapsto (i(s), \varphi (s, t)$.
Notice that $i$ is an involution of $S^{1}$. Moreover, since the
homology class of each circle $\set{ s  } \times S^{1}$ is an
eigenvector for the eigenvalue $1$ for the induced map by $f$ on
homology, $i$ must reverse the order and is therefore conjugate to
the map $s \mapsto -s$. In these new coordinates, $f$ is expressed
by $(s,t)\mapsto (-s, \phi_{s}(t))$ whereas $f^{2}$ is still the map
$(s,t ) \mapsto (s,t+2\alpha)$. Since $f$ commutes with $f^{2}$, so
does $\phi_{s}$ with the irrational rotation $t \mapsto t+2\alpha$.
Hence $\phi_{s}$ must also be a rotation and we can write
$\phi_{s}(t) = t+\alpha(s)$. The formula $f^{2}(s, t) = (s,
t+2\alpha)$ leads to $\alpha(s)+\alpha(-s) = 2\alpha$. In
particular, we have $2\alpha(0) = 2\alpha(1/2) = 2\alpha$ and hence
$\alpha(O) - \alpha(1/2) \in \set{0, 1/2 }$.

1) If $\alpha(0) = \alpha(l/2) = \alpha$, we let $\beta : S^{1}
\rightarrow S^{1}$ be the continuous map which is $0$ on $[0,1/2]$
and $\alpha(0) - \alpha(s)$ on $[-1/2, 0]$. The continuity of this
map at $1/2 = -1/2 \in S^{1}$ results from the fact that
$\alpha(0) = \alpha(l/2)$. An easy computation shows that for all
$s \in S^{1}$ we have:
\begin{equation}
\alpha(s) + \beta(s) - \beta(-s) = \alpha.
\end{equation}
We let then $B$ be the homeomorphism of the torus
defined by $B(s, t) = (s,t + \beta(s))$ and we can check that:
\begin{equation}
\left ( B^{-1}\circ f \circ  B \right )(s, t) = (-s, t + \alpha).
\end{equation}

2) If $\alpha(0) - \alpha(l/2) = 1/2$, we let $\beta: S^{1}
\rightarrow S^{1}$ be the continuous map which is s on $[0, 1/2]$
and $\alpha(0) -\alpha(s)$ on $[-1/2, 0]$. The continuity of this
map at $1/2 = -1/2 \in S^{1}$ results from the fact that
$\alpha(0) - \alpha(l/2) = 1/2$. An easy computation shows that
for all $s \in S^{1}$ we have:
\begin{equation}
\alpha(s) + \beta(s) - \beta(-s) = s + \alpha (0).
\end{equation}
Again, we let $B$ be the homeomorphism of the
torus defined by $B(s, t ) = (s,t + \beta(s))$ and we write $\alpha(0) = \alpha$.
An easy computation shows that:
\begin{equation}
\left ( B^{-1}\circ f \circ B \right ) (s, t) = (-s, t + s + \alpha).
\end{equation}
which completes the proof.
\end{proof}


\subsection{The case of the Klein bottle}\label{subsec:KleinBottle}

Let $\theta_{0}$ be the involution of the torus $\mathbb{T}^{2}$
defined by $(s,t ) \mapsto (-s, t + 1/2)$. Since this involution
is fixed point free, the quotient space
$\mathbb{T}^{2}/\theta_{0}$ is a closed surface and since
$\theta_{0}$ is orientation reversing, this surface must be the
Klein bottle $K$. The canonical map $\pi : \mathbb{T}^{2}
\rightarrow \mathbb{T}^{2}/\theta_{0}$ is the orientation covering
of $K$ and $\theta_{0}$ is the unique automorphism of this
covering.

For $\alpha \in (\mathbb{R}\setminus\mathbb{Q})/\mathbb{Z}$, we
note $\phi_{\alpha}$ and $\psi_{\alpha}$ the homeomorphisms of
$\mathbb{T}^{2}$ defined by $\phi_{\alpha}(s, t))= (s,t + \alpha)$
and $\psi_{\alpha}(s, t ) ) = ( s+ 1/2, t +\alpha)$. Since these
maps commute with $\theta_{0}$, they induces two homeomorphisms of
$K$ that we shall denote by $\Phi_{\alpha}$ and $\Psi_{\alpha}$,
respectively. These maps define two distinct families of regular
and non periodic homeomorphisms of $K$ and we note that:
\begin{enumerate}
    \item  For every irrational number $\alpha$, the closure of each $\Phi_{\alpha}$-orbit is a simple closed curve
which is the projection by $\pi$ of the two circles $\set{ s  }
\times S^{1}$ and $\set{-s } \times S^{1}$ of $\mathbb{T}^{2}$.
The rotation number of the restriction of $\Phi_{\alpha}$
restricted to this curve is equal to $\alpha$ if $s \neq -s$ (i.e.
$s \notin \set{ 0, 1/2  }$ or to $2\alpha$ if $s \in \set{0,1/2
}$.

    \item For every irrational number $\alpha$, the closure of each $\psi_{\alpha}$-orbit is the projection by $\pi$ of the
circles $\set{ s  } \times S^{1}$, $\set{ -s  } \times S^{1}$,
$\set{ s +1/2  } \times S^{1}$ and $\set{ -s + 1/2  } \times
S^{1}$ of $\mathbb{T}^{2}$. If $s \notin \set{ 1/4, -1/4 }$, these
circles are mapped into tow simple closed curve. If $s \in \set{
1/4, -1/4  }$, these circles are mapped into just one circle on
which the rotation number is $\alpha+ 1/2$.Note that
$\psi^{2}_{\alpha}=\Phi_{2\alpha}$.
\end{enumerate}

The aim of this paragraph is to establish the following
theorem:

\begin{thm} Let $f$ be a regular, non periodic
homeomorphism of the Klein bottle $K$ . Then there exist $\alpha \in (\mathbb{R}\setminus
\mathbb{Q})/\mathbb{Z}$ such that $f$ is topologically conjugate either to $\Phi_{\alpha}$ or to
$\Psi_{\alpha}$
\end{thm}

\begin{proof}
Recall first that $f$ has exactly two lifts on $\mathbb{T}^{2}$, one of which
$f_{+}$ preserves orientation while the other $f_{-}=\theta_{0}\circ f_{+}$ reverses
orientation. Note also that $f_{+}$ and $f_{-}$ commute with $\theta_{0}$.

Since $f_{-}$ and $\theta_{0}$ commute, so do the induced maps on
homology. According to the results of
section~\ref{subsec:TorusNegative}, the homology matrix of $f_{-}$
is therefore of type 1 and there are coordinates on the torus in
which $f_{-}$ can be written $(s,t ) \mapsto (-s, t+ \alpha)$. Since
$\theta_{0}$ commutes with $f_{-}$, its expression in these
coordinates must be $(s, t) \mapsto (i_{-}(s), t + \alpha(s))$ where
$i_{-}$ is an involution of $S^{1}$ which reverses the order. The
expression of $f_{+}$ follows from the relation $f_{+} = \theta_{0}
\circ f_{-}$ and we have $f_{+} (s, t) = (i_{+}(s), t + \alpha(s) +
\alpha)$ where $i_{+}(s) = -i_{-} (s)$ is an involution of $S^{1}$
which preserves the order.

It is a straightforward exercise to show that given two commuting
involutions of $S^{1}$, $i_{+}$ and $i_{-}$ such that $i_{+}$ preserves the order and
$i_{-}$ reverses the order, there is a homeomorphism of $S^{1}$ which
conjugate the pair $(i_{+}, i_{-})$ to either $(s \mapsto s, s \mapsto -s)$ or $( s \mapsto
s + 1/2,s \mapsto -s)$.

From this fact, we deduce that there are
coordinates on the torus in which we can write either:
\begin{description}
    \item[(1)] $f_{-}(s,t)=(-s,t+\alpha)$, $\theta_{0}(s,t)=(-s,t+\alpha (s))$, $f_{+}(s,t)=(s,t+\alpha (s) + \alpha)$
    \item[(2)] $f_{-}(s,t)=(-s,t+\alpha)$, $\theta_{0}(s,t)=(-s,t+\alpha (s))$, $f_{+}(s,t)=(s+1/2,t+\alpha (s) + \alpha)$
   \end{description}

Since $\theta_{0}$ is an involution, we have $\alpha(s) +\alpha(- s)= 0$. In particular $2\alpha(0) = 2\alpha(1/2)
= 0$ and we get $\alpha(0) = \alpha(l/2) = 1/2 \in S^{1}$ because $\theta_{0}$ is fixed point
free.

To achieve the proof of the theorem, we let $\beta : S^{1} \mapsto S^{1}$ be
the continuous map which is $0$ on $[0, 1/2]$ and $1/2 - \alpha(s)$ on $[-1/2,
0]$. The continuity of this map at $0$ and $1/2$ results from the fact
that $\alpha(0) = \alpha(1/2) = 1/2$. We verify then that, for all $s \in S^{1}$, we
have:
\begin{equation}
\alpha(s) + \beta(s) - \beta(-s) =1/2
\end{equation}
and we let $B$ be the homeomorphism of the torus defined by $B(s, t) = (s,t +\beta(s))$. In the first
case, we can check that:
\begin{itemize}
    \item $\left( B^{-1}\circ \theta_{0} \circ B \right) (s,t) = (-s, t + 1/2)$,
    \item $\left( B^{-1}\circ f_{+} \circ B \right) (s,t) = (s,t + \alpha + 1/2) = \Phi_{\alpha + 1/2}(s,t)$.
\end{itemize}
In the second case, we have $\alpha (s+ 1/2) = \alpha(-s)$ since $f_{+}$ commutes with
$\theta_{0}$ and we obtain:
\begin{itemize}
    \item $\left( B^{-1}\circ \theta_{0} \circ B \right) (s,t) = (-s, t + 1/2)$,
    \item $\left( B^{-1}\circ f_{+} \circ B \right) (s,t) = (s + 1/2,t + \alpha + 1/2) = \Psi_{\alpha + 1/2}(s,t)$.
\end{itemize}
which completes the proof.
\end{proof}


\subsection{The annulus and the M\"{o}bius
strip}\label{subsec:AnnulusMobius}

To obtain the complete classification of regular homeomorphisms of
the annulus $[-1, 1] \times S^{1}$ and the M\"{o}bius strip $[-1,1]
\times S^{1}/(s, t) )\mapsto (-s, t+1/2)$, it suffices to consider
their doubles and then to retain, among the list of regular
homeomorphisms of the torus and the Klein bottle, those which, keep
invariant an annulus or a M\"{o}bius strip. We state here the result
without any proof.

\begin{thm}
Let $f$ be a regular and non periodic homeomorphism of the
annulus. There exists $\alpha \in( \mathbb{R}\setminus \mathbb{Q})/\mathbb{Z}$ such that:
\begin{enumerate}
    \item $f$ $(s,t) = (s,t + \alpha )$ if $f$ is orientation-preserving,
    \item $f$ $(s,t) = (-s, t + \alpha )$
iff is orientation-reversing.
   \end{enumerate}
\end{thm}

\begin{thm}
Let $f$ be a regular and non periodic homeomorphism of the M\"{o}bius
strip. There exists $\alpha \in (\mathbb{R}\ \mathbb{Q})/\mathbb{Z}$
such that $f$ is conjugate to the map induced by the homeomorphism
$(s,t) \mapsto (s,t + \alpha)$ of the annulus.
\end{thm}

\bibliographystyle{amsplain}
\bibliography{regular}

\end{document}